\documentclass[10pt]{amsart}
\usepackage{amssymb,amsthm,amsmath}

\author{Benjamin Steinhurst}
\title[Dirichlet Forms on Laakso Spaces]{Dirichlet Forms on Laakso and Some Barlow-Evans Fractals of Arbitrary Dimension\footnote{Research supported in part by the National Science Foundation grant DMS-0505622.}}
\date{\today}

\newcommand{\C}[1]{\ensuremath{\mathcal{#1}}}
\newcommand{\B}[1]{\ensuremath{\mathbb{#1}}}

\newtheorem{theorem}{Theorem}[section]
\newtheorem{cor}{Corollary}[section]
\newtheorem{lemma}{Lemma}[section]
\newtheorem{prop}{Proposition}[section]

\newtheorem{defn}{Definition}[section]

\begin{document}

\begin{abstract}
In this paper we explore two constructions of the same family of metric measure spaces. The first construction was introduced by Laakso in 2000 where he used it as an example that Poincar\'e inequalities can hold on spaces of arbitrary Hausdorff dimension. This was proved using minimal generalized upper gradients. Following Cheeger's work these upper gradients can be used to define a Sobolev space. We show that this leads to a Dirichlet form. The second construction was introduced by Barlow and Evans in 2004 as a way of producing exotic spaces along with Markov processes from simpler spaces and processes. We show, for the correct base process in the Barlow Evans construction, that this Markov process corresponds to the Dirichlet form derived from the minimal generalized upper gradients. 


\end{abstract}

\maketitle

\section{Introduction}
\label{intro}

There is a sizable literature that considers fractal spaces and operators on them. A common simplification on the fractals to make the study more tractable is to assume that the fractals are finitely ramified, that is they can be disconnected by removing a finite number of points, \cite{Teplyaev2008,Bajorin2008}. A stronger but related simplification is to consider fractals that are post-critically finite as done in \cite{Kigami2001, BarlowPerkins1988, Kusuoka1989, Lindstrom1990}. There has also been interest in post-critically infinite fractals in \cite{BarlowBass1989, BarlowBass1999, Bajorin2008}. Often these spaces can be very poor in paths between points leading to problems in conducting analysis on them \cite{Heinonen2007, Semmes1996}. The main obstacle this presents is that it prevents the use of capacity and curve modulus arguments to obtain Poincar\'e inequalities and other related objects. 

In \cite{Heinonen2007} there is an excellent survey of the kinds of analysis which can be done on spaces which are not smooth in a classical sense but which do still have a ``large supply'' of rectifiable curves connecting any two points. One of the notable results which can come from having enough curves in a space is a $(1,1)-$Poincar\'e inequality. Unfortunately, many fractals do not have this ample supply of curves, for example the Sierpinski Gasket. Laakso, in \cite{Laakso2000}, gave a construction of an one-parameter family of metric measure spaces which have sufficient rectifiable curves to support a Poincar\'e inequality with the advantage that the (continuous) parameter indexing the family of spaces is the Hausdorff dimension. The dimension of the space is determined by the number and dimension of the Cantor sets as well as the sequence $\{ j_i \}$ used in the construction that we review in Section\ \ref{LConst}. Moreover, a countable subfamily of these Laakso spaces are self-similar fractals. Laakso's original construction is an elegant one but not well suited to studying the properties of operators on these spaces. Examples of spaces other with Ahlfors regularity and probabilistic information such as escape time estimates, as opposed to analytic information, are discussed in \cite{Barlow2004}.  

The spaces that Laakso constructed have enough rectifiable curves to allow for the kind of analysis in \cite{Semmes1996} which uses the capacity of sets as a central tool. We will define a Dirichlet form on each of the Laakso spaces that is derived from the minimal generalized upper gradients of suitable functions. Barlow and Evans \cite{BarlowEvans2004} have constructed Markov processes that evolve on what they call ``vermiculated spaces,''  and state that Laakso's spaces can be constructed as vermiculated spaces. Starting with the Dirichlet form associated to the minimal generalized upper gradient we identify the Markov process to which it corresponds. 

In \cite{BarlowEvans2004},  there are proofs of the existence of Markov processes on Barlow-Evans spaces using a construction which we show can give Laakso spaces, although it can generate a much wider variety of spaces as well. This leads to a natural question: whether these Markov processes are symmetric with respect to a reasonable measure on the space? If there are symmetric processes the next question would be to which Dirichlet forms do they correspond? The connection between Dirichlet Forms and Markov processes is well known and we refer to the exposition from \cite{FOT1995} for the general theory.  Following up on a comment in \cite{BarlowEvans2004} we offer a proof that Laakso's spaces can also be constructed as projective limits of quantum graphs. We surmise that Barlow and Evans knew this result but did not include it in their paper. This perspective will be used in Section\ \ref{SharedProc} to prove the final theorem of the paper.  In \cite{Teplyaev2008} a similar use of quantum graphs to estimate a Dirichlet form is explored for finitely ramified fractals based on \cite{Kigami1993,Kigami2003}; our situation is much more complicated though. 

Both Dirichlet forms and symmetric Markov processes are associated to unbounded, self-adjoint operators. Once we have proved that we can realize Laakso spaces as projective limits of quantum graphs we will be able to show that the operator associated to the minimal generalized upper gradient Dirichlet form is also the limit of the operators on the sequence of approximating quantum graphs that are associated to a particular Markov process taken through the Barlow-Evans construction. In this way we will produce a symmetric Markov process and a Dirichlet form on any Laakso space which correspond to the same operator, hence are associated themselves. We analyze the spectra of such operators in \cite{RomeoSteinhurst2009} and \cite{BegueEtAl2009}.

We begin by reviewing the basic theory of Dirichlet forms and Markov processes on general spaces in Section \ref{DFormsMProc}. Then in Section \ref{LConst}\ we give in detail Laakso's original construction from \cite{Laakso2000}. In Sections \ref{SF-Section} and  \ref{DFUP-Section} we define a space of functions and then describe explicitly a Dirichlet form and minimal upper gradients on the fractal. The construction offered by Barlow and Evans in \cite{BarlowEvans2004} is presented in Section \ref{BEConst}\ and sufficient conditions for the existence of Markov processes on the vermiculated spaces are given in Section\ \ref{ProcBE}. Then in Section \ref{SharedProc}\ we link this Dirichlet form to a specific Barlow-Evans Markov process. 

{\bf Acknowledgments:} This paper would not have come into being without the guidance and support of Alexander Teplyaev. I also thank Piotr Haj\l asz, Luke Rogers, and Robert Strichartz for their useful questions and comments. Comments from anonymous referees have proved invaluable in improving this paper.

\section{Dirichlet Forms and Markov Processes}
\label{DFormsMProc}
 In this section we briefly recall some basic facts about the relation between these two approaches and establish notation. The reader can find more details in \cite{FOT1995,RogersWilliams2000}. While our approach is mainly analytic in much of the previous literature, including \cite{BarlowEvans2004}, probabilistic approaches have been used.

There is a deep connection between Dirichlet Forms, which are on the face primarily analytic objects, and Markov Processes, which are very much probabilistic objects. This connection is classical and has been explored by many authors including Fukushima, Oshima, and Takeda in \cite{FOT1995}. We begin our discussion by recalling basic definitions and stating without proof a theorem that gives the conditions necessary for the correspondence. We assume that all Hilbert spaces mentioned in this paper are real $L^{2}$ spaces on the relevant space. Throughout this section we assume a regular measure space $(X,\mu)$.

\begin{theorem}\label{SFSG-Cor}
There is a one-to-one correspondence between closed symmetric bilinear forms on a Hilbert space and non-positive definite self-adjoint operators on the Hilbert space. The correspondence is given by:
$$\left\{ \begin{array}{l} Dom(\C{E}) = Dom(\sqrt{-A})\\ \C{E}(u,v) = (\sqrt{-A}u,\sqrt{-A}v).\end{array}  \right. $$
\end{theorem}

The correspondence is between \C{E} and $-A$ where for any $u \in Dom(A)$, $\C{E}(u,u) = (u,-Au)$. And $Dom(A) \subset Dom(\sqrt{-A})$ is a dense, proper subset. See \cite{FOT1995} Thm 1.3.1 for the proof or \cite{Rudin1973} is another standard reference.

This $-A$ is an operator on the underlying Hilbert space which can be viewed as the generator of a semi-group via $exp(tA)$ or alternatively as the generator of a resolvent via $(\alpha - A)^{-1}$ where these expressions are given meaning by a spectral resolution and the functional calculus for self-adjoint operators. Naturally this induces correspondences between closed symmetric forms, operators, semi-groups of operators, and resolvents. 

\begin{defn} A \emph{Dirichlet Form}, $(\C{E}, Dom(\C{E}))$, is a closed bilinear symmetric form on an $L^{2}(X,\mu) = L^{2} = H$ space such that if $u \in Dom(\C{E}) \subset H$ then $(u \vee 0) \wedge 1 \in Dom(\C{E})$ and $\C{E}((u \vee 0) \wedge 1,(u \vee 0) \wedge 1) \le \C{E}(u,u)$. This type of contraction of $u$ is called a unit contraction.
\end{defn}

In \cite{FOT1995} instead of the $(u \vee 0) \wedge 1$ being again in the domain the authors use a differentiable function $\phi(x)$ being in the domain where $\phi(x) = x$ for $x \in [0,1]$, $\phi(x) \in [-\epsilon, 1+\epsilon],\ 0 \le \phi'(x) \le 1,\ \forall x \in \B{R}$. And contract $u(x)$ by composition with $\phi(x)$. This type of contraction is a \emph{normal} contraction. However these conditions are equivalent when the form is closed. 

The adjective \emph{Markovian} applies to operators, semi-groups of operators, and symmetric forms. Due to the connection between semi-groups and symmetric forms the usages correspond to each other but before we state that correspondence we set out what those properties are in each case:
\begin{itemize}
	\item Bounded Linear Operator: An operator, $S$, is \emph{Markovian} if for all $0 \le u \le 1$ it is the case that $0\le Su \le 1$ where the inequalities hold almost everywhere.
	\item Semi-group: A semi-group of bounded linear operators, $\{ T_t, t \ge 0\}$, is \emph{Markovian} if for all $t \ge 0$ the operator $T_t$ is Markovian.  
	\item Symmetric Form: A symmetric form, $D$, is \emph{Markovian} if for all $\epsilon >0$ there is a non-decreasing function $\phi_{\epsilon}(t)$ such that $\phi_{\epsilon}(t)=t$ if $t \in [0,1]$, $-\epsilon < \phi_{\epsilon}(t) < 1+\epsilon$, and $\phi_{\epsilon}(t') - \phi_{\epsilon}(t) \le t'-t$, and $u\in Dom(D) \Rightarrow \phi_{\epsilon}(u) \in Dom(D).$ If $Dom(D)$ is closed this is equivalent to the unit contraction $(u \vee 0) \wedge 1 \in Dom(D)$ and $\C{E}((u \vee 0) \wedge 1, (u \vee 0) \wedge 1) \le \C{E}(u,u)$ and is defined above.
\end{itemize}

Notice that a Markovian symmetric form has all the properties of a Dirichlet form except being closed. 

\begin{theorem}\label{DFMP-Cor}
Let \C{E} be a closed symmetric form on $L^{2}(X,m)$. Let $\{ T_t, t>0\}$ and $\{G_{\alpha},\alpha>0\}$ be the strongly continuous semigroup and the strongly continuous resolvent on $L^{2}(X,m)$ which are associated with $\C{E}$. Then the following are equivalent:
\begin{enumerate}
	\item $T_t$ is Markovian for each $t>0$.
	\item $\alpha G_{\alpha}$ is Markovian for each $\alpha>0$.
	\item \C{E} is Markovian, i.e. a Dirichlet form.
	\item For any $u \in Dom(\C{E})$, $(u \vee 0) \wedge 1 \in Dom(D)$ and $\C{E}((u \vee 0) \wedge 1, (u \vee 0) \wedge 1) \le \C{E}(u,u)$. This is referred to as the unit contraction ``operating'' on \C{E}.
	\item For any $u \in Dom(\C{E})$, $\phi_{\epsilon}(u) \in Dom(\C{E})$ and $\C{E}(\phi_{\epsilon}(u), \phi_{\epsilon}(u)) \le \C{E}(u,u)$. This is referred to as the normal contraction ``operating'' on \C{E}.
\end{enumerate}
\end{theorem}

See \cite{FOT1995} Thm 1.4.1 for the proof.

This theorem states that the use of the word \emph{Markovian} in these different settings is an appropriate use of terminology. At the end of the next group of definitions and theorems these contexts will be connected to stochastic processes in which the word \emph{Markovian} was first used. 

We define the basic probabilistic objects and notation that we will need to be able to state which processes the Dirichlet forms will correspond. Denote by $\Omega$ a sample space, $\C{F}$ a $\sigma-$field on $\Omega$, $X_t$ a process which is adapted to the filtration $\C{F}_t \subset \C{F}$, $P^{x}$ is the law of $X_t$ when $X_0=x$. Denote by $S$ the state space of $X_t$ with Borel field, $\C{B}$. Adjoin a point, $\Delta$, to $S$ to serve as a \emph{cemetery} point. Let $S_{\Delta} = S \cup \{\Delta\}$ and $\C{B}_{\Delta} = \C{B} \cup \{ B \cup \Delta : B \in \C{B}\}$.  Later in the paper we will use the Laakso fractals and approximations to them as state spaces. 

\begin{defn} A quintuplet $(\Omega, \C{F}, \C{F}_t,  \{X_t\}, \{P^{x}\} : t \in [0,\infty], x \in S_{\Delta})$ is a \emph{Markov process} if the following conditions hold:
\begin{enumerate}
	\item The quintuplet is a progressively measurable stochastic process with $t$ as the time parameter and $(S_{\Delta}, \C{B}_{\Delta})$ as its state space. 
	\item There exists an admissible filtration $\{\C{M}_t\}_{t \ge 0}$ which has the property that for each $x \in S$, $$P^{x}(X_{s+t} \in E|\C{M}_t) = P^{X_t}(X_s \in E)\ a.s.$$ For any $s,t \ge 0$ and $E \in \C{B}$. 
	\item $P^{x}(X_t \in E)$ is \C{B}-measurable as a function of $x$ for all $t \ge 0$ and $E \in \C{B}$ and $P^{x}(X_0 = x)=1$.
	\item $P^{\Delta}(X_t = \Delta) =1$ for all $t \ge 0$.
\end{enumerate}
\end{defn}

Adjoining the cemetery point compactifies the state space and assures that the associated symmetric form is conservative. However, in our case all of the state spaces will already be compact and we can take the cemetery as an inaccessible state and still have conservative symmetric forms. 

To each Markov process, $X_t$, associate the transition function $p_t$ where $p_t(x,A) =  P^{x}(X_t \in A)$ where $A$ is a Borel subset of $S$ which acts on functions by $p_t u(x) = \int u(y)p_t(x,dy)$. If only a single probability measure is given on the sample space $\Omega$, then one can use a similar definition $p_t(x,A) =\B{P}(X_t(\omega) \in A | X_0(\omega) = x)$ where $\B{P}$ is the given probability measure. For each $t > 0$ $p_t(x,A)$ is a kernel, and a \emph{Markovian} kernel if $p_sp_t = p_{s+t}$ and $0 \le p_t(x,A) \le 1$. Then $\{p_t\}$, integration against which is a symmetric operator, generates a semi-group of symmetric integral operators on $L^{2}$ for each $t >0$, called $T_t$. We will need strongly continuous semigroups for the correspondence to Dirichlet forms, so to ensure that $T_t$ is strongly continuous at zero we have the following criterion:

\begin{lemma} If $p_t(x,A)$ is a symmetric Markovian transition function and $T_t$ the associated semi-group of operators. Then $T_t$ is strongly continuous at zero if $\lim_{t \downarrow 0} p_tu(x) = u(x)$ for $u$ that are continuous with compact support in $S$. 
\end{lemma}

This next theorem gives the next piece of the correspondence. But first we need another definition.

\begin{defn}
A \emph{Hunt process} is a Markov process which almost surely has right continuous and left quasi-continuous sample paths. See \cite{FOT1995} for more on quasi-continuity. A \emph{Diffusion} is a Markov process that almost surely has continuous sample paths. A Hunt process or diffusion is symmetric if its infinitesimal generator is a symmetric operator, or equivalently the associated heat kernel is symmetric in the spacial coordinates.
\end{defn}

\begin{theorem}\label{SGMP-cor} There is a one to one, up to equivalence, correspondence between symmetric Markovian transition semi-groups and symmetric Hunt processes. 
\end{theorem}

\begin{proof}
This is a combination of Theorems 7.2.5 and 4.2.7  of \cite{FOT1995}.
\end{proof}

The correspondence between a Hunt process and a Dirichlet form is through the semi-group generated by the process and the associated infinitesimal generator. This generator is the operator defined by the Dirichlet form, the $-A$ in the notation used in the definition of Dirichlet form above. Since we are often interested in looking at processes with continuous sample paths we note that continuity of sample paths translates along the correspondence to the Dirichlet form having the local property.

\begin{defn} A Dirichlet form, \C{E} is \emph{regular} if the compactly supported continuous functions in $Dom(\C{E})$ are dense in $Dom(\C{E})$ under the $\|u\| = \C{E}(u,u) + (u,u)$ norm, and $Dom(\C{E})$ is sup norm dense in the space of compactly supported continuous functions on the underlying space.

A Dirichlet form, $\C{E}$, possesses the \emph{local property} if for any $u,v \in Dom(\C{E})$ which have compact, disjoint support $\C{E}(u,v) = 0$.  
\end{defn}

We now have all the pieces to be able to state the final and most specific correspondence that we will mention in this section. 

\begin{theorem}\label{DFMP-Cor2}
The following two conditions are equivalent to each other for a regular Dirichlet form \C{E} on $L^{2}(X,\mu)$:
\begin{enumerate}
	\item \C{E} possesses the local property.
	\item There exists a $\mu$-symmetric diffusion process on $(X, \C{B}(X))$ whose Dirichlet form is the given one, \C{E}.
\end{enumerate}
\end{theorem}

This is Theorem 7.2.2 from \cite{FOT1995}.

Another property of Markov processes will become important later in the paper so we give the definition of Feller processes here. Let $$C_{\infty}(X) = \{ f \in C(X) : \forall \epsilon >0,\ \exists K\ compact, |f(x)| < \epsilon, \forall x \in X \setminus K\}.$$ This is the space of functions vanishing at infinity. When $X$ is itself compact $C_{\infty}(X) = C(X)$.

\begin{defn}\label{def:Feller}\cite{FOT1995} A Markov process is \emph{Feller} if the associated resolvent $G_{\lambda}$ has the following property:
$$G_{\lambda} C_{\infty}(X) \subset C_{\infty}(X)$$
That is that the resolvent maps continuous functions vanishing at infinity into continuous functions vanishing at infinity. for all $\lambda > 0$.
\end{defn}

There are variations on the definition of Feller processes in the literature for example in \cite{RogersWilliams2000} it is the semi-group $P_t$ that is considered and not the resolvent. The Hille-Yosida theorem which gives the relation between resolvents and semi-groups shows that the two approaches yield the same results.

\section{Laakso Construction}
\label{LConst}
This construction was first presented in \cite{Laakso2000} as a way to provide examples of metric-measure spaces with nice analytic properties e.g. a Poincar\'e inequality, of any arbitrary Hausdorff dimension greater than one. The original treatment made no mention of any probabilistic structures associated with the constructed space, though minimal upper gradients were shown to exist. For more on minimal upper gradients see \cite{Cheeger1999}. All of these spaces will have a cell structure, and for a countable collection of $Q$ the cell structure will be self-similar, making these self-similar spaces fractals.

We first mention a few facts about Cantor sets. The standard Cantor set can be constructed with two affine contraction mappings. One, $\psi_1$, mapping the interval $[0,1]$ to $[0,\frac13]$ and the other, $\psi_2$, mapping $[0,1]$ to $[\frac23,1]$. Then the Cantor set can be defined as the unique non-empty compact subset of $\mathbb{R}$, $K$, such that $K = \psi_1(K) \cup \psi_2(K)$. The Cantor set has Hausdorff dimension $\ln(2)/\ln(3)$ where the two is the number of contraction mappings and the one third the contraction factor, see \cite{Falconer1990} for more about the dimension of self-similar sets. One can change the Hausdorff dimension by altering the contraction factor to be anything in $(0,\frac12)$. The cell structure of Cantor sets is defined below. The properties of the cell structure, associated contraction mapping, and exactly calculable Hausdorff dimension extend to products of Cantor sets. 

\begin{defn}\label{cantorcell}
The zero level cell is $K$, the entire Cantor set. If $\psi_1$ and $\psi_2$ are the contraction mappings that define the Cantor set $K$ by the relation $K = \psi_1(K) \cup \psi_2(K)$ then $K$ is the zero-level cell, $\psi_i(K)$ is a first level cell, $\psi_i(\psi_j(K))$ is a second level cell, and so on. A \emph{cell} of a Cantor set is a cell of any level. The \emph{cell structure} of a Cantor set is the set of all the cells of every level. 
\end{defn}

For a given dimension, $Q>1$, we begin with two spaces. The first is a Euclidean space, $I = [0,1]$. The second is a product of cantor sets, $K^{k}$ where each $K$ has Hausdorff dimension $\frac{Q-1}{k}$ so that the product has dimension $Q-1$. Consider the product space $I \times K^{k}$, where the measure is the product of the Lebesgue measure on $I$ and the product Bernoulli measure on $K$. Note that $I \times K^{k}$ has total measure one. The fractal, $L$, will be the quotient space of $I \times K^{k}$ by an equivalence relation where the identifications will be made on a null set so that there will be a natural, induced measure $\mu$ on $L$ that is Borel regular. 

To be able to find where the identifications will be made we need a number derived from the desired dimension of $K$. Let $t \in (0,1)$ such that $\ln(2)/\ln(1/t) = \frac{Q-1}{k}$ where $k$ is chosen large enough so that $\frac{Q-1}{k} \in (0,1)$. This gives a $t$ to be used as the contraction factor in the iterative construction of the Cantor set, the fraction of the length of an interval at the $m$th step that the intervals at the $m+1$st level are. This gives a natural decomposition of $K = tK \cup (tK+1-t)$. When we take the product $I \times K^{k}$ it will have dimension $Q$. It is necessary to have a way of describing the location of a point in the Cantor sets with an ``address.'' Call $K_0 := tK$ and $K_1:= tK+1-t$ then $K_{00}$ is the left part of the left part of $K$ i.e. $t^{2}K$. This naming scheme can be continued and associates to each point $x \in K$ an address $a = a_1a_2a_3 \ldots$ so that $x = K_a$. Finite addresses indicate subsets of $K$ and can be concatenated to produce the addresses of still smaller subsets. If $a$ is a finite address let $|a|$ be its length. This scheme for labeling the points of a self-similar space is heavily used in \cite{Kigami2001}.

For the given $t$ there exists an integer $j$ such that $\frac{1}{j+1} < t \le \frac{1}{j}$. Then there is a sequence $j_i \in \{j, j+1\}$ such that 
\begin{equation}\label{eqn:jt}
	\frac{j}{j+1} \prod_{i=1}^{m} j_i^{-1} \le t^{m} \le \frac{j+1}{j} \prod_{i=1}^{m} j_i^{-1}.
\end{equation}
Now define a function $w$ which will pinpoint exactly where each level of identifications will occur.

\begin{defn}\label{def:w}
For $l \ge 1$, define the function
\begin{equation}\label{eqn:w}
	w(m_1, \cdots, m_l) = \sum_{i=1}^{l}m_i \prod_{h=1}^{i} j_h^{-1}
\end{equation}
Where  $0 \le m_i < j_i$ for $i <l$ and when $i=l$ $0 < m_l <j_l.$ The values of $w(m_1, \cdots, m_l)$ give the locations of the $l'th$ level wormholes in the $I$ coordinate.
\end{defn} 

The condition on $m_l$ forces the wormholes to not stack up on each other by forbidding $l$-level wormholes from being located over any lower level wormholes. Suppose that there are $k$ Cantor sets used in constructing a particular Laakso space, then we consider the set of points in $I$ with coordinates taken from the values of $w(m_1, \cdots, m_l)$ and let $(x_1,x_2, \ldots, x_{k+1})$ be a point with first coordinate in $I$ and the rest of the coordinates in $K^{k}_{a0}$ where $a$ is an address with length $k-1$. We identify $(x_1,x_2,\ldots, x_{k+1}) \in I \times K_{a0}^{k}$ with $(y_1,y_2,\ldots, y_{k+1}) \in I \times K_{a1}^{k}$ if and only if the length of $a$ is $l-1$, $x_1 = y_1$ is a value of $w(m_1, \cdots, m_l)$, and $y_i = x_i + t^{l-1}(1-t).$ We make these identifications iteratively for all $l$. The points at which these identifications are made are known as \emph{wormholes}. 

\begin{defn}\label{def:iota} Denote the identification map sending $I \times K^{k}$ to $L$ by $\iota$.
\end{defn}

The space $L$ has, by construction, a cell structure already in the $K^{k}$ as each one of these Cantor sets has the normal cell structure. If the $j_i$ are also periodic, then there is self-similarity in the $I$ direction as well. Say that the $j_i$ have period $p$, then let a cell of $L$ be the image under the identification map of the set $[m/r, (m+1)/r] \times K_{a_i}^{k}$ where $m = 0 ,1, \cdots r$ where $r = \prod_{l=1}^{p} j_l$ and the $a_i$ are addresses of length $p$ and there is one (potentially different) address for each copy of $K$ used. Any function defined on $L$ can be defined on $I \times K^{k}$ as a pullback by the identification map. So a function $f:L \rightarrow \mathbb{R}$ can alternatively  be worked with as $\hat{f}(x,w) = f \circ \iota:I \times K^{k} \rightarrow \mathbb{R}, x \in I, w \in K^{k}$ as well. 

A simple approach to showing a space metrizable is to construct a metric. The most natural metric on this space is a geodesic metric where the distance between two points is the infimum of lengths of all rectifiable paths connecting the two points. The existence of rectifiable curves connecting any two points in $L$, which implies that the space is connected and that the geodesic metric is well defined, is shown in \cite{Laakso2000}. 

Laakso's construction gives an easy to use measure, namely the product measure on $I \times K^{k}$ carried down by the identification map. This measure is also the $Q$-Hausdorff measure on $L$. We now summarize the basic properties of $L$ before moving onto defining function spaces. 

\begin{theorem}
The space $L$ is a connected metric measure space which is Alfors regular of dimension $Q$.
\end{theorem}

This is the central result of \cite{Laakso2000}. 

It is worth taking some care in understanding how the geodesic metric behaves on $L$. The length of a rectifiable path comes entirely from the distance that it travels  in the $I$ direction since traversing a wormhole to move from one copy of $I$ to another costs no length. One can then use the arc length parameterization of a path to induce a measure on the image of that path. These measures are the one dimensional Lebesgue-Stieltjes measures associated to the rectifiable paths. Call these measures $dm$, but keep in mind that they are dependent on the specific path over which the integral is taken.

\begin{defn}\label{MGUG} On a metric measure space, $(X,|\cdot |)$, a \emph{minimal generalized upper gradient} of a function $u$ is a non-negative function $p_u$ with the following property: $$|u(x) - u(y)| \le \int_{\gamma} p_u\ dm$$ 
For any pair of points $x,y \in X$ and rectifiable curve $\gamma(t) \subset X$ such that $\gamma(0) = x$ and $\gamma(1) = y$ and any other function with this same property is almost everywhere greater than or equal to $p_u$ and the measure $dm$ is the measure induced by $\gamma$.
\end{defn}

It is a simple matter to note that the function $p = \infty$ is a generalized upper gradient. Thus a generalized upper gradient exists for any function. We follow Cheeger \cite{Cheeger1999} in viewing the set of functions in $L^{p}(L)$ that have a generalized upper gradient also in $L^{p}(L)$ as a Sobolev space, $H^{1,p}$. If $p>1$ then there exists a unique minimal generalized upper gradient. 
A more complete overview of abstract Sobolev spaces is at the end of the next section. 
It is more convenient to be able to speak of only one upper gradient, this is fine so long as $p>1$.

\begin{theorem} \cite{Cheeger1999} 
For $1<p<\infty$ if $f \in H^{1,p}$ there exists a minimal generalized upper gradient which is unique up to modification on sets of measure zero.
\end{theorem} 

The intuition behind this theorem is that for $p>1$ $L^{p}$ is a convex space, so minimizing sequences of generalized upper gradients actually have a unique limit point. The rest of the proof is checking that the limit is again a generalized upper gradient.

Note that in a Euclidean space for differentiable functions the minimal generalized upper gradient is the norm of the usual gradient, $p_u = |\nabla u|$, so in a sense $p_u$ plays the same role as the absolute value of a more general first derivative. With this generalized minimal upper gradient we have, from \cite{Laakso2000}, a weak $(1,1)-$Poincar\'e inequality: $$\int_B |u-u_B|\ d\mu \le C(diam(B))\left(\int_{CB} p_u\ d\mu \right).$$ Here $B \subset L$ is a ball, $\mu$ is the measure on $L$, and $C$ is a constant.

\section{A Space of Smooth Functions}
\label{SF-Section}
In this section we define a space of functions, $\mathcal{G}$, on $L$ which will serve as a core for the Dirichlet form that is defined in the next section. We then prove that for these functions the minimal generalized upper gradient is easy to describe. For notational simplicity we assume that $k=1$, that only one Cantor set is being used in the construction. Finally we define a Sobolev space based on minimal generalized upper gradients.

\begin{defn}\label{def:hatnotation}
For a function $f \in C(L)$ denote by $\hat{f}(x,w):I \times K \rightarrow \mathbb{R}$ the pulled back function $f \circ \iota$ (c.f. Definition \ref{def:iota}).
\end{defn}

\begin{defn}\label{def:calG}
If $f \in C(L)$ and for $\hat{f}(x,w)$ there exists an $n \ge 0$ such that when $K$ is decomposed into cells of depth $n$ and these two conditions are met:
\begin{itemize}
	\item for a fixed $x \in I$ the function $\hat{f}(x,\cdot)$ is constant on each cell of $K$ (see Definition \ref{cantorcell});
	\item for a fixed $w\in K$ the function $\hat{f}(\cdot,w)$ is continuously differentiable between wormhole locations of depth $n$ or less with finite limits at the wormhole locations;
\end{itemize}
then we say that $f \in \mathcal{G}_n$. Let $\mathcal{G} = \bigcup_{n=0}^{\infty} \mathcal{G}_n$.
\end{defn}

When we pull back to a function $f \circ \iota = \hat{f}(x,w)$ on $I \times K$ the infinitely many identifications are already accounted for by having started with precisely those functions in $f \in C(L)$. The main point of this definition is to be able to analyze $f \in C(L)$ in terms of it's ``directional'' behavior which doesn't become well defined until $f$ is pulled back to $\hat{f}(x,w) = f \circ \iota$. Also in this definition when we define $\mathcal{G}_n$ and force the functions to be constant on each $n$th level cell for a given $x \in I$ we are in effect treating the product of an interval between two wormhole locations crossed with cell of $K$ as a single line segment making $I \times K$ look like a quantum graph. Increasing $n$ then increases the complexity of this graph allowing for more functions on $L$ that are included. This intuition will be revisited in Section \ref{SharedProc}.

\begin{lemma}\label{StoneWeierstrass}
The space $\mathcal{G}$ is dense in the continuous functions on $L$ in the supremum norm.
\end{lemma}

\begin{proof}
To use the Stone-Weierstrass Theorem we need to show that the algebra $\mathcal{G}$ separates points and contains the constant functions. Constant functions are all elements of $\mathcal{G}_0 \subset \mathcal{G}$. Let $p,q \in L$ be distinct points. Then they either have different coordinates in the $I$ direction or they don't. If they do then $\hat{f}(x,w) = x \in \iota^{*}\mathcal{G}_0$ and will separate the points $p$ and $q$. The space $\iota^{*} \mathcal{G}_0$ consists of all of the pull backs of functions in $\mathcal{G}_0$ to functions on $I \times K$, so there exists a function $f \in \mathcal{G}$ such that $\hat{f} = x$ and it will separate $p$ and $q$. If $p$ and $q$ have the same coordinate in the $I$ direction, say $x_0 \in I$, then they must have different coordinates in the $K$ direction which can be distinguished by cells of some finite level, say $p =\iota( p(x_0,w_1))$ and $q=\iota(q(x_0,w_2))$ where $w_1,w_2 \in K$ are in different $n$th level cells of $K$, call them $K_1$ and $K_2$. To construct a separating function $\hat{f}(x,w)$ in this case if $p$ and $q$ are in different $n$th level cells and not at a wormhole of level $n$ or lower then there are wormholes with locations $y < x_0 < z$ such that they are the closest to $x_0$. Let
$$\hat{f}(x,w) = \left\{ \begin{array}{rl} (x-y)(x-z) & x \in (y,z), w \in K_1\\ -(x-y)(x-z)& x \in (y,z), w\in K_2\\ 0 & otherwise \end{array}\right. .$$ 
Then $\hat{f}(x,w)$ is defined on finite level cells and is piecewise defined from differentiable functions $\hat{f}(x,w) \in \iota^{*} \mathcal{G}$ and $f(p) = -f(q) \neq 0$. If $x_0$ is the location of an $n$th level wormhole simply use the same process with $n+1$ level cells of $K$. That this function is the pull back of a well defined function on $L$ holds by checking that the first $n$ levels of identifications in the construction of $L$ are respected and that any lower level of identification are as well.

The regions in $L$ defined by $\iota([y,z] \times K_i)$ are actually $n$th level cells of $L$ because the number of wormholes at each level can be chosen randomly $(\{j_i\}$ need have no pattern) $L$ is not necessarily a self-similar fractal so there aren't analogues to the $\psi_i$ in Definition \ref{cantorcell} to be used in defining a cell structure. We use a notion of cell structure based on Definition 2.1 in \cite{Teplyaev2008} that does not rely on self-similarity. In this notion the cells are a family of subsets for each scale of $L$, $\{F_{\alpha}\}_{\alpha \in A}$ along with a family of boundaries $\{B_{\alpha}\}_{\alpha \in A}$, where $F_{\alpha} \cap F_{\alpha'} = B_{\alpha} \cap B_{\alpha'}$. This condition states that the intersection of two cells is the intersection of their boundary. The situation in \cite{Teplyaev2008} is one where these boundaries are finite sets of vertices, but each boundary of a cell in a Laakso space is a Cantor set. To see why the cells $\iota([y,z] \times K_i), i=1,2$ have disjoint interiors in $L$ it becomes necessary to know how a path from a point in $\iota(I \times K_1)$ could reach a point in $\iota(I \times K_2)$. When we defined the identification maps $n$th level cells could only be connected by $n$th level and lower (i.e. $n-1$ level) wormholes so if no such wormholes are in the interior of cells then the cells can at most share their boundaries which is no problem in defining our function since it is zero on the boundary of the two cells so when these sets are mapped back into $L$ their interiors remain disjoint.
\end{proof}

\begin{theorem}\label{GradientLemma1} For $f \in \mathcal{G}_n \subset \mathcal{G}$ and $q$ not a wormhole i.e. $\iota^{-1}(q) = (x,w)$ and $x \neq w(m_1,m_2, \ldots, m_k)$ for any $k \le n$ (see Definition \ref{def:w} and following for the definition of this function), $p_f(q) = \left| \frac{\partial}{\partial x} \hat{f}(\iota^{-1}(q)) \right|$ where $x$ is the $I$ coordinate in $I \times K$ and $q \in L$ for $\mu$-a.e. $q \in L$. \end{theorem}

The set of wormholes forms a set of measure zero and are ignored since minimal generalized upper gradients are only defined almost everywhere. As a short hand we denote $\left| \frac{\partial}{\partial x} \hat{f}(\iota^{-1}(q)) \right|$ as  $|\frac{\partial}{\partial x}f |$. 

\begin{proof}
First we show that $\left| \frac{\partial}{\partial x}\right|$ is a generalized upper gradient then we show that it is the minimal one. Given the boundedness assured i nthe definition of $\mathcal{G}$, this upper gradient is also integrable and square integrable. Now take two points $x,y \in L$ and a rectifiable path connecting them, $\gamma$. When $\gamma$ is pulled back to $\tilde{\gamma}$ on $I \times K$ there is ambiguity at each wormhole that $\gamma$ goings through so make the choices that make the lifted $\tilde{\gamma}$ right continuous and have left limits in the time parameter. Because $f \in \mathcal{G}$ it is associated to a decomposition of $K$ into cells of some finite level. Then even if $\tilde{\gamma}$ is completely disconnected it must have some length in the $I$ direction in each cell that is passes through. This is because the only way a wormhole can provide a path out of an $n$th level cell is for the wormhole to be at most of depth $n-1$ which are evenly spaced. Let $x=z_0$ is in one of the cells of $K$, let $z_1$ be the point in $L$ when $\tilde{\gamma}$ first leaves this cell, $z_2$ the point in $L$ when $\tilde{\gamma}$ first leaves that that cell, and so on. Since $\gamma$ is a rectifiable path it has finite length and $\gamma$ will only transit finitely many of these cell at most countably many times, with possible repeats. So let $z_m$ be the $m^{th}$ crossing from one cell to another and $z_{\infty} = y$. It may happen that the path only moves from one cell to another finitely many times, in that case the modification is obvious. 

In each $n$th level cell of $K$ the requirement that $\hat{f}(x,w)$ be constant across the cell for a given $x \in I$ means that in each $n$th level cell of $K$ $\hat{f}(x,w)$ is a piecewise differentiable function in $x$. This means that along the path $\tilde{\gamma}$, as it passes through a cell, standard calculus methods can be used to determine an upper gradient in that cell, which will be the usual $\left| \frac{\partial}{\partial x}f\right|$. Since this can be done for each of the countably many transts that $\tilde{\gamma}$ makes through the cells and in fact for any $\tilde{\gamma}$ that we may have chosen we use on all cells the generalized upper gradient $\left| \frac{\partial}{\partial x}f\right|$ on all of $L$

To show minimality we proceed by contradiction. Suppose that there is another generalized upper gradient, $p_f$ which is less than  $\left| \frac{\partial}{\partial x}f\right|$ on a set of positive measure, $A$. Then there is a subset $A'$ of $A$ with positive measure such that $\frac{\partial}{\partial x}f$ is of one sign and this subset contains open sets by the piece wise continuity of the derivatives of functions in $\mathcal{G}$. There is a subset $A''$ of $A'$ that is contained in a single cell of $K^{k}$ crossed with some subinterval of $I$, without loss of generality assume that $\frac{\partial}{\partial x}f$ is positive. Then for $p_f$ to be less than or equal to $\frac{\partial}{\partial x}f$ on a set of positive measure would imply that on one-dimensional intervals that the absolute value of the first derivative is not the minimal generalized upper gradient which is a contradiction.
\end{proof}

\begin{defn}\label{H12} The \emph{Sobolev space} $H^{1,2} \subset L^{2}(L)$ is defined to be $$H^{1,2} = \{ u \in L^{2}(L)| \exists p_u,  p_u \in L^{2}(L)\}$$ Where $p_u$ is the minimal upper gradient of $u$. $H^{1,2}$ is given the graph norm: $$\|u\| = \left(\int u^{2}\right)^{1/2} + \left(\int p_u^{2}\right)^{1/2}.$$
\end{defn}

This definition is from \cite{Cheeger1999} where the following lemma is proved (Theorems 2.7 and 2.10).

\begin{lemma}\label{lem:cheegerprop}
The space $H^{1,2}$ is a complete Banach space and the minimal generalized gradients are unique up to modification on a set of measure zero.
\end{lemma}

\begin{lemma}\label{SobolevG} The Sobolev space $H^{1,2}$ contains the closure of $\C{G}$ under the graph norm. That is the set $$\bar{\C{G}} \subset H^{1,2}  = \{u \in L^{2}(L)| \exists p_u, p_u \in L^{2}(L)\}.$$
\end{lemma}

\begin{proof} Since the Sobolev space $H^{1,2}$ is complete the only thing that needs to be checked is that $\mathcal{G} \subset H^{1,2}$. Since all elements of \C{G} have bounded derivatives (see Definition \ref{def:calG}) on a finite measure space we see that $\int u^{2} < \infty$ and $\int p_u^{2} < \infty$ so $\|u \| < \infty$ hence $\mathcal{G} \subset H^{1,2}$.
\end{proof}

\begin{theorem}\label{thm:diffu}
For any $u \in \overline{\mathcal{G}}$ the map $u \mapsto \frac{\partial}{\partial  x} u \in L^{2}(L)$ is well-defined.
\end{theorem}

\begin{proof}
For any function $u \in \mathcal{G}$ the object $\frac{\partial}{\partial x}u$ exists by passing to the pull back, $\hat{u}(x,w)$ where the definition of $\mathcal{G}$ assure us that $\frac{\partial}{\partial x}u$ exists a.e.. Now let $u_n \in \mathcal{G}$ such that $u_n \rightarrow u$ in $H^{1,2}$. This convergence implies that $u_n \rightarrow u$ in $L^{2}$ and $\frac{\partial}{\partial x}u_n$ is Cauchy in $L^{2}$. This means that $\frac{\partial}{\partial x}u_n$ converge to an unique element in $L^{2}(L)$ which we will call $\frac{\partial}{\partial x}u$. In \cite{Cheeger1999} it is shown that $p_u = \lim_{n \rightarrow \infty} \left| \frac{\partial}{\partial x}u_n \right|$ is actually the minimal generalized upper gradient of $u$ so $\left| \frac{\partial}{\partial x} u \right| = \lim_{n \rightarrow \infty} \left| \frac{\partial}{\partial x}u_n \right| = p_u$. Thus the relationship between $p_u$ and $\frac{\partial}{\partial x}u$ that was observed for $u \in \C{G}$ extends to $\overline{\C{G}}$.
\end{proof}

It seems reasonable that the inclusion $H^{1,2} \subset \overline{\mathcal{G}}$ holds as well, but we do not need it for this paper. It rests on the consideration of whether $\bar{\C{G}}$ is dense in the Lipschitz functions, \cite{Hajlasz2003}. This is still an unresolved issue. We continue by defining our Dirichlet form on \C{G} and then take the closure of \C{G} as the domain.

\section{A Dirichlet Form and Upper Gradients}
\label{DFUP-Section}

In this section we show how to use generalized minimal upper gradients to produce a Dirichlet form. We use the space of functions \C{G}, see Definition \ref{def:calG}, on which the generalized minimal upper gradients can explicitly be computed. It would be a natural choice to define a Dirichlet form $\C{E}(u,u) = \int_L p_u^{2}\ d\mu$ and then use polarization to extend to a bilinear form that would look like $\C{E}(u,v) = \int_L p_u p_v\ d\mu.$ This can't be done because $u \mapsto p_u$ is not a linear operator, so setting $p_u = \sqrt{-A}u$ is not a  viable definition. Recalling the notation from Section \ref{DFormsMProc}, we have another option at our disposal. We have the map $u \mapsto \frac{\partial}{\partial x}u$ that can be used instead and since this map is linear the extension to a bilinear form will hold. 

\begin{lemma}\label{LaaDF} Let $u \in \bar{\C{G}} = Dom(\C{E})$, let $$\C{E}(u,u) = \int_L p_u^{2}\ d\mu = \int_L \left( \frac{\partial}{\partial x} u \right)^{2}\ d\mu.$$ Then $(\C{E},Dom(\C{E}))$ is a Dirichlet form. 
\end{lemma}

{\bf Remark:}\label{rem:A} If $u \in \C{G}$ is such that $\hat{u}(x,w)$ is piecewise twice-differentiable in the $I$ direction with derivatives vanishing on the boundary and have directional derivatives summing to zero at the wormholes. Integration by parts with suitable boundary conditions imposed indicates that we can also express the Dirichlet form as $$\C{E}(u,u) = -\int_L u \frac{\partial^{2}u}{\partial x^{2}}\ d\mu.$$ The domain of $-\frac{\partial^{2}}{\partial x^{2}} = -A$ is a strict, but dense, subset of $\bar{\C{G}} = Dom(\C{E})$.


\begin{proof}
By defining $\C{E}$ through polarization it can be seen that \C{E} is a symmetric, quadratic form on $H^{1,2}(L,\mu)$.  That \C{E} is closed is ensured by the general construction of the ambient Sobolev space in  \cite{Cheeger1999} which addressed Dirichlet forms constructed from minimal generalized upper gradients. This leaves the Markov property to check. Since $\C{E}$ is a closed form it suffices to check the Markov property on $\mathcal{G}$. We shall use Theorem \ref{DFMP-Cor} to use a normal contraction instead of the unit contraction as in our definition. Since the normal contraction $\phi_{\epsilon}$ is a differentiable function if $u \in \mathcal{G}$ then $\phi_{\epsilon} \circ u \in \mathcal{G}$ with $p_u$ being a generalized upper gradient for $\phi_{\epsilon} \circ u$ so that $\C{E}(\phi_{\epsilon} \circ u,\phi_{\epsilon} \circ u) \le \C{E}(u,u)$.  It is worth noting that $p_{\phi_{\epsilon} \circ u}$ is in general point-wise less than $p_u$.
\end{proof}

By this lemma we see that the Dirichlet form $(\C{E}, Dom(\C{E}))$ is generated by the self-adjoint operator $-A$ whose domain is a dense subspace of $Dom(\C{E})$. At present it is unclear whether $Dom(\C{E}) = H^{1,2}$ since this rests on the density of $\C{G} \subset H^{1,2}$ (see the comments after Lemma \ref{SobolevG}). However many Dirichlet forms have smaller domains than the ambient Sobolev space so this is not an unusual situation. One could view $Dom(\sqrt{-A})$ as a sort of first order Sobolev space and $Dom(A)$ as a second order Sobolev space. Care must be taken when using this analogy to remember that these spaces are embedded in, but not equal to $H^{1,2}$. The operator, $-A$, comes back in to consideration at the end of the paper. Now we show two properties of \C{E}, locality and regularity. 

\begin{theorem}\label{DFSummary} The symmetric form \C{E} is a local, regular Dirichlet whose domain is contained in the  Sobolev space $H^{1,2}$ and is the closure of the function space \C{G} under the graph norm associated to the operator $\sqrt{-A}$. For the function spaces this implies that $$\overline{\mathcal{G}} = Dom(\sqrt{-A}) = Dom(\mathcal{E}).$$
\end{theorem}

\begin{proof}
The only two things left to check are locality and regularity. Locality holds for $u \in \mathcal{G}$ as an immediate consequence of the definition of $\mathcal{G}$ and Theorems \ref{GradientLemma1} and \ref{thm:diffu}.  Because $\C{E}$ is closed this can be extended to the entire domain of the form. Regularity is a consequence of Lemma \ref{StoneWeierstrass}\ and the definition of $H^{1,2}$. 
\end{proof}

\begin{cor} There is a non-trivial Markov process with continuous sample paths on the fractal $L$.
\end{cor}

\begin{proof}
Non-triviality follows from $\C{E}(f,f) \neq0$ when $\hat{f}(x,w) = x^{2}$. Continuity of sample paths is from Theorem\ \ref{DFMP-Cor2}\ which requires both the locality and regularity that we have established. 
\end{proof}

We end this section with an overview of the various definitions of Sobolev spaces with mention of various equivalencies. Much of this discussion is taken from \cite{Heinonen2007}. We  begin with a quick summary of what has already been done in this paper. 

What we've done here has been to find a replacement for the norm of the gradient in building up Sobolev spaces on Laakso spaces.  This brings the classical notions of Sobolev spaces which may be stated in terms of the Laplacians associated to Dirichlet forms defining the Sobolev space. In the remark after Lemma \ref{LaaDF} we stated what the Laplacian and its domain are. With these two objects in hand we recall the various notions of Sobolev spaces on Euclidean domains.

Historically Sobolev spaces began with function spaces over domains in $\mathbb{R}^{n}$. Such as $W^{1,2}(\mathbb{R})$ which is the space of square integrable functions with square integrable first derivatives. This forces the members of $W^{1,2}$ to have a desired amount of smoothness. Higher derivatives could be required or $p$-integrability instead of square integrability to get spaces $W^{1,p}$ this only changes to exponent in the integrability condition and nothing else in the definition. In the light of distribution theory we might want to ease the smoothness requirement and require that the distributional derivatives be integrable instead of classical derivatives, these spaces are also known as $H^{1,p}$. These spaces coincide with $W^{1,p}$ when the boundary of the domains is suitably smooth, so over all of $\mathbb{R}^{n}$ or on disks they are the same \cite{Heinonen2001}. These definitions still restrict us to spaces locally reminiscent of $\mathbb{R}^{n}$ to be able to talk about ``derivatives'.' 

Then Haj\l asz \cite{Hajlasz1996, FranchiHK1999} extended the concept into arbitrary metric measure spaces. His definition was for a space $M^{1,p}$ which consisted of all functions $u$ for which there existed another function $g$ such that for all $x,y $ in the space 
$$|u(x)-u(y))| \le d(x,y)(g(x)+g(y)).$$ 
With $g$ acting as a sort of maximal function since there need not even exist such a $g$ for an arbitrary $u$ and the norm is $\|u\| = \|u\|_p + \inf \|g\|_p$ where the infimum is taken over all $g$ with the required property. In Euclidean spaces $M^{1,p} = W^{1,p}$ \cite{Heinonen2001}. But when upper gradients are introduced it allows another definition of a Sobolev space in a metric- measure space. Shanmugalingam \cite{Shanmugalingam2000} introduced Newtonian Spaces where a function $u$ is in $N^{1,p}$ if there exists some $p_u$ such that 
$$|u(x)-u(y)| \le \int_{\gamma} p_u\ dm$$ 
Where this must hold for some function $p_u \in L^{p}$ and for all but a capacity zero set of paths $\gamma$ connecting $x$ and $y$. It is known that if the space supports a $(1,q)$-Poincar\'e inequality then Sahnmugalingam and Haj\l asz's Sobolev spaces coincide \cite{Hajlasz2003}. We use in this paper Cheeger's \cite{Cheeger1999} version of this type of space which also relies on upper gradients, which is the definition already given of the Sobolev space above.

\section{Barlow-Evans Construction}
\label{BEConst}

In \cite{BarlowEvans2004}, Barlow and Evans commented that their construction can also produce Laakso's spaces. They do not, however, prove this statement. This fact is very useful because it lends itself to providing alternative proofs for the existence of Dirichlet forms and Markov processes on Laakso spaces. Towards this end we describe their construction, prove that the Laakso spaces can be constructed this way, and show that there are many Dirichlet forms on Laakso's spaces. Barlow and Evans' construction is based on Evans' previous work with Sowers in \cite{EvansSowers2003}.

To construct a vermiculated space, $L$, one needs three ingredients. The first is a state space, $F_0$, for the base Markov process. We'll take $F_0 = [0,1]$ to construct Laakso spaces. The second is a family of sets, $G_n$, which at each step of the construction will index the possible \emph{alternate universes} or copies of $F_{n-1}$ that the process could evolve in, these sets are taken to be $\{0,1\}$ to construct Laakso spaces with dimension less than two. The last ingredient is another family of sets, $B_n$, which indicate where the identifications or ``wormholes'' between the  $\#G_n$ copies of $F_{n-1}$ are made. It is the sequence $\{B_n\}_{n=1}^{\infty}$ that will determine the dimension of $L$.

We begin with $F_0 = F$, and the sequences $\{G_n\}_{n=1}^{\infty}$, and $\{B_n\}_{n=1}^{\infty}$. The construction is inductive. Define $E_1 = F_0 \times G_1$, $\hat{E}_1 = F_0$, and $A_1 = B_1 \times G_1 \subset E_1$. Note that $B_1 \subset F_0$. The next two functions are defined so as to perform the identifications that will create the next approximation to $L$, namely $F_1$. Define $\psi_1:E_1 \rightarrow \hat{E}_1 = F_0$ by $\psi_1(y,z) = y.$ Let $\tilde{E}_1 = (E_1 \setminus A_1) \cup \psi_1(A_1)$ with the topology induced by the function $$\pi_1(y,z) = \left\{ \begin{array}{cc} (y,z), & if\ y \in F_0 \setminus B_1,\\ y, & if\ y \in B_1;\end{array} \right.$$ Let $F_1:= \tilde{E}_1$. There is also a continuous surjection $\phi_1:F_1 \rightarrow F_0$ given by $$ \begin{array}{cc} \phi_1 (y,z) = y, & (y,z) \in E_1 \setminus A_1 = (F_0 \setminus B_1) \times G_1\\ \phi_1(y) = y, & y \in \psi_1(A_1) = B_1. \end{array}$$ This construction can be repeated by using $F_1, G_2, B_2$ to produce $F_2$ and so on. The set $\{F_n\}_{n=0}^{\infty}$ along with the surjections $\{\phi_{n,m}\}_{n,m=1}^{\infty}:F _m \rightarrow F_n$ form a projective system whose inverse limit is the space $\lim_{\leftarrow} F_i \subset \prod_{n=1}^{\infty} F_n$. Due to basic facts from \cite{HockingYoung1988} about projective limit spaces $\lim_{\leftarrow} F_i$ is compact and Hausdorff since all of the $F_n$ are compact and Hausdorff.

\begin{figure}[tbp]\centering
\begin{picture}(200,120)
\put(90,115){$\prod F_i$}
\put(85,70){$\lim_{\leftarrow} F_i$}
\put(100,110){\vector(0,-1){30}}
\put(90,110){\vector(-1,-2){45}}
\put(110,110){\vector(1,-2){45}}
\put(90,65){\vector(-1,-1){45}}
\put(110,65){\vector(1,-1){45}}
\put(35,10){$F_n$}
\put(160,10){$F_m$}
\put(50,10){\vector(1,0){105}}
\put(90,00){$\phi_{n,m}$}
\put(70,30){$\Phi_n$}
\put(120,30){$\Phi_m$}
\end{picture}
\caption{Summary of the Projective System, $n > m \ge 0$}
\label{ProjSystem}
\end{figure}
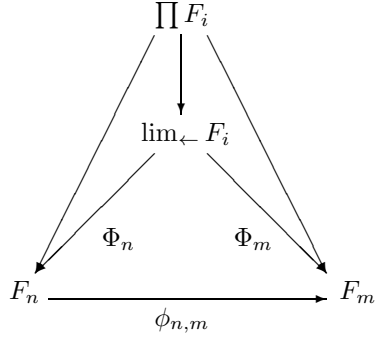

The projective system is summarized in Figure\ \ref{ProjSystem}. In the proof of the following Lemma we describe each of the maps explicitly and then use the Universal Property of Projective Limits to show that the map $\eta$ in Figure\ \ref{UnivProp} is an isometry between the Laakso space and the inverse limit space of the $F_n$. In our opinion these considerations are best understood in conjunction with the example instead of in abstract terms in the preceding discussion. 

The claim of the following Lemma is one of the primary goals of this section, proving that the Barlow-Evans construction can be used to construct Laakso spaces and thus is a more general construction. One of the benefits of using the Barlow-Evans construction is that the existence of Markov processes is also ensured. After this Lemma we will show that there is not one, but many of these Markov processes on $L$. We believe that Barlow and Evans were aware of a proof of this fact but have not published it.

\begin{lemma} There exists a homeomorphism, $\eta$ between any Laakso fractal and a vermiculated space. 
\end{lemma}

\begin{proof}
First we show that there exists a continuous surjection, $\eta$, from a given Laakso space onto a particular Barlow-Evans space which is constructed in the following paragraphs. Then we show that $\eta$ is also injective with continuous inverse. 

This proof is specialized for Laakso spaces with dimension between 1 and 2 to simplify notation. However, for higher dimensions use products of the $G_n$ that we define. Take $F_0 = [0,1]$ and $G_n = G = \{0,1\}$. For a given $t \in (0,\frac12)$ one can find a sequence of $j_m \in \{j,j+1\}$ where $j\le t^{-1} < j+1$. This sequence should be chosen such that 
$$\frac{j}{j+1}\prod_{i=1}^{m}j_i^{-1} \le t^{m} \le \frac{j+1}{j}\prod_{i=1}^{m}j_i^{-1}.$$ 
Note that this is the same sequence of integers that was chosen in the Laakso construction above, see Equation \ref{eqn:jt}. Let $B_n$ consist of points  of the form 
$$w(m_1,m_2,\ldots , m_n) = \sum_{i=1}^{n} m_i \prod_{h=1}^{i}j_h^{-1}$$ 
Where $0 \le m_i \le j_i$ with the additional proviso that $m_n >0$, and $g_n$ any point in $G_n$, this is the same function that gave the location of the wormholes in the Laakso construction in Definition \ref{def:w} and following. These choices will put wormholes at the same locations in the Barlow Evans construction. This makes $A _n = B_n \times G_n \subset F_{n-1} \times G_n$ as needed. If each $F_{n-1} \times G_n = \hat{E}_n$ is taken to lie in the unit square in $\B{R}^{2}$ with lower left corner at the origin then then horizontal slices are approaching a Cantor set of the necessary dimension and vertical slices are copies of the unit interval \cite[Section 2-14]{HockingYoung1988}.

Inductively construct the spaces $F_i$ as described above. These spaces come with maps $\phi_{i+1,i}:F_{i+1} \rightarrow F_i$. Let $\phi_{i,j} = \phi_{j+1,j} \circ \cdots \circ \phi_{i-1,i}:F_{i} \rightarrow F_j$ for $i >j$. Now consider the space $\prod_{i=0}^{\infty} F_i$, there are projection maps which we will call $\Phi_n:\prod_{i=0}^{\infty} F_i \rightarrow F_n$ for all $n \ge 0$. The projective limit space $\lim_{\leftarrow}F_i$ will actually be a subspace of $\prod_{i=0}^{\infty} F_i$ that we can explicitly define. Define $\lim_{\leftarrow} F_i$ to be all elements $\{x_i\}_{i=0}^{\infty}  \in \prod_{i=0}^{\infty} F_i$ such that $\phi_{i,j}(x_i) = x_j$ for all $i>j \ge 0$ \cite[Page 91]{HockingYoung1988}. We can then restrict the maps $\Phi_n$ to $\lim_{\leftarrow} F_i$ since it is a subspace of $\prod_{i=0}^{\infty} F_i$ and we will call the restrictions $\Phi_n$ as well leaving it to context to make it clear which space they project from. It is important to note how $\Phi_n$ and $\phi_{i,j}$ interact since they are all projection operators we have that $\Phi_j = \phi_{i,j} \circ \Phi_i$ for $i>j$. Projective limits are a very general concept that is even treated in Category theory. The most pertinent property of projective limit systems is the \emph{universal property.} This property is a statement that in a certain sense the projective limit is minimal. Minimality in this sense means that if another topological space $L$ has maps $\tilde{\Phi}_n:L \rightarrow F_n$ for all $n \ge 0$ such that $\tilde{\Phi}_j = \phi_{i,j} \circ \tilde{\Phi}_i$ for all $i>j$ that there is an induced a continuous surjection $\eta:L \rightarrow \lim_{\leftarrow} F_i$ that factors $\tilde{\Phi}_i$ as $\Phi_i \circ \eta.$ The diagram in Figure \ref{UnivProp}. Moreover, this diagram commutes. To show that we can take $L$ to be a given Laakso space we need to construct $\tilde{\Phi}_i$ such that  $\tilde{\Phi}_j = \phi_{i,j} \circ \tilde{\Phi}_i$ and then we will know that $\eta$ is a continuous surjection from a given Laakso space onto the Barlow Evans space constructed to have the same wormholes as the Laakso space. 

Let $L$ be a Laakso space, then define $\tilde{\Phi}_i:L \rightarrow F_i$ to be given by $\iota_i \circ (id,\pi_i) \circ \iota^{-1}$ where $\iota$ is the identification map sending $I \times K$ in the Laakso construction, $(id, \pi_i):I \times K \rightarrow I \times \C{K}$, \C{K} is the collection of endpoints of the depth $i$ cells of the Cantor set $K$ (a finite set), there is only one copy of $K$ since we have restricted ourselves to spaces with dimension less than two. And $\iota_i$ is the identification map that only identifies the wormholes of level $i$ or less. Refer to Section \ref{LConst} for the original discussion of the Laakso space's construction. The composition of these maps is continuous and surjects onto a quantum graph that can be identified with $F_i$. Since $\pi_i \circ \pi_j = \pi_i$ we will have that $\tilde{\Phi}_j = \phi_{i,j} \circ \tilde{\Phi}_i$ for $i>j\ge 0$. By the universal property of projective limits then map $\eta:L \rightarrow \lim_{\leftarrow} F_i$ exists and is a continuous surjection.

\begin{figure}[tbp]\centering
\begin{picture}(200,150)
\put(150,140){$\prod F_i$}
\put(150,80){$\lim_{\leftarrow} F_i$}
\put(30,80){$L$}
\put(160,130){\vector(0,-2){40}}
\put(45,83){\vector(1,0){100}}
\put(95,87){$\eta$}
\put(30,10){$F_n$}
\put(155,10){$F_m$}
\put(35,75){\vector(0,-1){50}}
\put(35,75){\vector(2,-1){110}}
\put(20,45){$\tilde{\Phi}_n$}
\put(65,60){$\tilde{\Phi}_m$}
\put(45,13){\vector(1,0){100}}
\put(160,75){\vector(0,-1){50}}
\put(160,75){\vector(-2,-1){110}}
\put(110,60){$\Phi_n$}
\put(170,45){$\Phi_m$}
\put(90,00){$\phi_{n,m}$}
\end{picture}
\caption{Use of the Universal Property}
\label{UnivProp}
\end{figure}
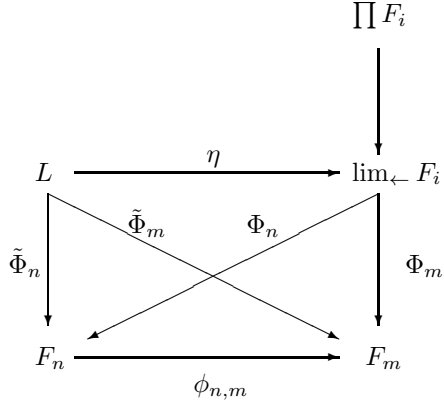

It remains to be proved that $\eta$ is injective and has continuous inverse. Both families $\Phi_i$ and $\tilde{\Phi}_i$ separate points. For $\Phi_i$ this is because if they didn't separate two points the construction of the projective limit would have made them the same point. For the $\tilde{\Phi}_i$ it is because any two points in $L$ can eventually be distinguished by cells in the Cantor set of some finite level or by $I$-coordinate. Suppose that $\eta$ is not injective then there exists distinct $p,q \in L$ such that $\eta(p) = \eta(q)$. Thus $\Phi_i(\eta(p)) = \Phi_i(\eta(q))$ for all $i \ge 0$. But $\tilde{\Phi}_i(p) \neq \tilde{\Phi}_i(q)$ for some $i$ since $p \neq q$. Since the diagram in Figure \ref{UnivProp} is a commutative diagram $\Phi_i\circ \eta = \tilde{\Phi}_i$ and we have a contradiction to the commutativity of the diagram. Thus $\eta$ is bijective. Since $\eta$ is a continuous bijection from a compact Hausdorff space into a Hausdorff space, it is a homeomorphism.
\end{proof}

We have shown that $\eta(L) = \lim_{\leftarrow}F_i$ and that $\eta$ is a homeomorphism so for the rest of the paper we will simply say that $L = \lim_{\leftarrow} F_i$ and identify the function spaces as well.

Since we are interested in processes on the limit space we need to also consider a projective system of measures as well. Recall that $\Phi_n:L \rightarrow F_n$ is a projection from the limit space to the $nth$ approximating space, then $\Phi_n^{*}: B(F_n) \rightarrow B(L)$ maps the functions spaces by composition i.e. $\Phi_n^{*}(f) = f \circ \Phi_n^{*}$. If we use $\Phi_n^{*}$ to map indicator functions we can use $\Phi_n^{*}$ to map sets from the finite approximation spaces to the limit space. In our example on each $F_n$ there is a measure, $\mu_n$, that is a weighted one-dimensional Lebesgue measure on the Quantum graph with total mass one. Alternatively it can be viewed as the measure induced in the quotient space $F_n$ by Lebesgue measure on $F_{n-1} \times G_n$. To be a projective system of measures the collection $\{\mu_n\}$ must be compatible
$$\mu_{n+1}(\phi_{n+1,n}^{*}A) = \mu_n(A)$$
For $A \in B(F_n)$ and $\phi_{n+1,n}^{*}:B(F_n) \rightarrow B(F_{n+1})$ defined the same way as $\Phi^{*}_n$. The $\mu_n$ have bounded total mass so by \cite[Prop 8, III.50]{Bourbaki2004} there is a unique limit measure such that $\mu_{\infty}(\Phi_n^{*} U) = \mu_n(U)$ if $U$ is a measurable subset of $F_n$. The concern will be if this measure $\mu_{\infty}$ can be given in concrete terms adapted to our situation. That is represent it in a way such that it can be worked with. Here we show it to be the same measure as obtained from the Laakso construction. 

\begin{lemma}  Let $\mu$ be the measure obtained in the Laakso construction and $\mu_{\infty}$ the measure obtained from the Barlow-Evans constructiion. Then $\mu = \mu_{\infty}$.
\end{lemma}

\begin{proof}
Since we know that the spaces $\lim_{\leftarrow}F_i$ and $L$ are topologically the same and we have a subbasis for the topology in both, which generates the $\sigma-$algebra on which the measures are defined. Since the measures are finite, as long as they agree on the algebra generated by the basis elements the measures will agree on all measurable sets. We take as subbasis elements $(r,s) \times K_{a_1} \times \cdots \times K{a_n}$ where $r,s$ are \emph{not} wormholes and the $a_i$ are finite length addresses. The intersection of two elements of this subbasis is again an element of the subbasis. Call the measure on the Laakso construction $\mu$ to distinguish it from $\mu_{\infty}$. Then $\mu((r,s) \times K_{a_1} \times \cdots \times K{a_n}) = |r-s| 2^{-|a_1|} \cdots 2^{-|a_n|}$, where $|a_i|$ is the length of the address $a_i$. If the maximum length of the $a_i$ is $M$ then $(r,s) \times K_{a_1} \times \cdots \times K{a_n}$ is the image under $\Phi_M^{*}$ of some rectangle-like set in $F_M$ which has $\mu_M$ measure $|r-s| \times 2^{-|a_i|} \cdots 2^{-|a_n|}$. Since $\mu_M$ is the product measure with identifications on a set of measure zero it agrees with $\mu$. Since these sets generate the Borel $\sigma-$algebra both $\mu$ and $\mu_{\infty}$ are extensions of the same finite pre-measure and so are equal. Thus the the map $\eta$ from the proof of the previous theorem is a measure preserving homeomorphism.
\end{proof}

Before moving onto considering random processes on these two spaces we take advantage of the measure preserving isometry $\eta$. It is a recapitulation of the preceding results. 

{\bf Remark:} The spaces $L$ and $\lim_{\leftarrow} F_i$ are identified through the map $\eta$. The Sobolev space $H^{1,2}(L)$ is naturally identified with a function space on $\lim_{\leftarrow} F_i$ via composition with the map $\eta^{-1}:\lim_{\leftarrow} F_i \mapsto L$ which is also called $H^{1,2}$. Similarly with any function space, such as $L^{2}$, $Dom(A)$, or $Dom(\C{E})$, on either $L$ and $\lim_{\leftarrow} F_i$ are identified. Since $\eta$ is an isometry between $L$ and $\lim_{\leftarrow} F_i$ then for $f:L \rightarrow \mathbb{R}$ we have $f \circ \eta : \lim_{\leftarrow} F_i \rightarrow \mathbb{R}$ and for $g: \lim_{\leftarrow} F_i \rightarrow \mathbb{R}$ we have $g \circ \eta^{-1} : L \rightarrow \mathbb{R}$. Because $\eta$ is itself a continuous bijection pre-composing with $\eta$ or $\eta^{-1}$ a function space on either $L$ or $\lim_{\leftarrow} F_i$ can be viewed as a function space on the other. 

The Sobolev space $H^{1,2}$ is defined in Definition \ref{H12}. In Lemma \ref{LaaDF}, the spaces $Dom(\mathcal{E})$ and $Dom(A)$ are defined as the domains of the Dirichlet Form on the Laakso construction and the domain of the associated Laplacian. Note that this is not yet enough to show that any of the function spaces, other than the continuous functions, defined via the Laakso construction or the Barlow Evans construction coincide this is addressed in the remaining sections.

\section{Processes on Barlow-Evans Spaces}
\label{ProcBE}

In \cite{BarlowEvans2004}, Barlow and Evans present not only a construction of state spaces using projective limits but also sufficient conditions on a base Markov process on $F_0$ so that a Markov process on the limit space can be constructed. They show this process to be a Hunt process. We maintain the notation from the previous section concerning the names of sets involved in the Barlow-Evans construction, but from now on we'll only consider the process on $F_0$ to be reflected Brownian motion on the unit interval. 
 
\textbf{Assumption:} Write \C{C} for the collection consisting of the empty set and finite unions of sets drawn from $B_1, B_2, \ldots$. Assume that for each $C \in \C{C}$ that the resolvent of the process $X_t^{0}$ stopped on hitting $C$ maps $C(F_0)$ into itself.

We use this assumption in the context where $X^{0}_t$ is the Markov process on $F_0$ that we wish to extend. This assumption is much stronger than saying that $X^{0}_t$ is a Feller process and will allow us to show that $X^{n}_t$ is Feller as well. We note that standard Brownian motion on a line fits the assumption but Brownian motion on the plane does not \cite{BarlowEvans2004} if the sets in \C{C} are finite sets of singletons.

\begin{prop}\label{1D-BR} One-dimensional reflected Brownian motion on the unit interval  satisfies the assumption with the $B_i$ being finite point subsets of the unit interval.
\end{prop}

\begin{proof}
Given any finite set of points in the unit interval, $B$, and a Brownian motion starting at any point and stopped at $B$ the Brownian motion will behave, including its resolvent, like Brownian motion on an interval of finite, and possibly zero, length where the endpoints stop the process. Since Brownian motion has continuous sample paths it cannot escape from between which ever two points of $B$ it started between. Thus as long as Brownian motion stopped at end points has a resolvent that maps continuous functions to continuous functions this assumption will be satisfied.  

The resolvent map as defined by $f \mapsto (\alpha - \Delta)^{-1}f = g$ with Dirichlet boundary conditions to describe the absorbing boundaries of the process on this interval is an ODE which has a differentiable solution. So $(\alpha- \Delta)^{-1}C(F_0) \subset C(F_0)$.
\end{proof}

The process of constructing a sequence of Markov processes on the space $F_n$ is a repeated use of the method set forth in \cite{EvansSowers2003} whereby the process $X^{n+1}_t$ on $F_{n+1}$ is constructed from $X^{n}_t$ by extending the resolvents $U_{\alpha}^{n}$ associated to $X^{n}_t$ to be resolvents $U_{\alpha}^{n+1}$ on $F_{n+1}$. These resolvents $U^{n+1}_{\alpha}$ are then associated to a Markov process which is called $X^{n+1}_t$ which evolves on $F_{n+1}$. The limiting process which gives $U^{\infty}_{\alpha}$, the resolvent for the limit process, is described by Barlow and Evans in \cite{BarlowEvans2004}. We are going to use Theorem \ref{DFMP-Cor2}\ to link the limit process, $X^{\infty}_t$, to a Dirichlet form, $\tilde{\C{E}}$, that can be compared to the Dirichlet form, \C{E}, from Section \ref{DFUP-Section}. The hypotheses of Theorem \ref{DFMP-Cor2} must be checked, the first of which is symmetry of the process.

\begin{lemma}\label{lemma:symmetry} The Markov process on $L$ built from a Markov process $X_t$ on $F_0$ is symmetric with respect to the measure $\mu_{\infty}$ if $X^{0}_t$ is symmetric on $F_0$ with respect to Lebesgue measure on the unit interval. 
\end{lemma}

\begin{proof}
It follows by construction that $\bigcup_{n} \Phi_n^{*} C(F_n)$ is dense in $C(L)$ \cite{BarlowEvans2004}, where $C(F_n)$ are the continuous functions on $F_n$. So to talk about the symmetry of $U_{\alpha}^{\infty}$ it is sufficient to consider only functions from $\bigcup_{n} \Phi_n^{*} C(F_n)$. Then if $f,g \in \bigcup_{n} \Phi_n^{*} C(F_n)$ we have $f = \Phi_N^{*} \tilde{f}$ and $g = \Phi_N^{*} \phi_{M,N}^{*}\tilde{g}' = \Phi_N^{*}\tilde{g}$ where $\tilde{f}, \tilde{g} \in C(F_N)$. The value of $N$ is simply indicating at which level of approximation both $f$ and $g$ are describable without loss of information. Let $U_{\alpha}^{n}$ be the resolvent associated to the process $X^{n}_t$ on $F_n$ and $U^{\infty}_{\alpha}$ be the resolvent associated to the process $X_t$ on $L$. The relation that defines $U^{\infty}_{\alpha}$ on $\bigcup_{n} \Phi_n^{*} C(F_n)$ is 
\begin{eqnarray}\label{eqn:dynkin}
	U^{\infty}_{\alpha}\Phi^{*}_n f &=& \Phi^{*}_n U^{n}_{\alpha} f \hspace{1cm} \forall f \in C(F_n),\ \forall n \ge 0.
\end{eqnarray}
This relation defines $U^{\infty}_{\alpha}$ on $\Phi^{*}C(F_n)$ for every $n$ which since $\bigcup_{n} \Phi_n^{*} C(F_n)$ is dense in $C(L)$ $U^{\infty}_{\alpha}$ can be extended by continuity to all of $C(L)$. This relationship between $U^{\infty}_{\alpha}$ and $U^{n}_{\alpha}$ is called the Dynkin Intertwining relationship. Then by the Dynkin Intertwining relationship that holds for these resolvents we have:
\begin{eqnarray*}
	(f,U_{\alpha}^{\infty}g)_L &=& (U_{\alpha}^{\infty}f,g)_L\ iff\\
	(\Phi_N^{*}\tilde{f}, U_{\alpha}^{\infty} \Phi_N^{*}  \tilde{g})_L &=& (U_{\alpha}^{\infty} \Phi_N^{*}\tilde{f}, \Phi_N^{*} \tilde{g})_L\ iff\\
	(\Phi_N^{*}\tilde{f}, \Phi_N^{*}U_{\alpha}^{N}  \tilde{g})_L &=& (\Phi_N^{*}U_{\alpha}^{N} \tilde{f}, \Phi_N^{*} \tilde{g})_L\ iff\\
	(\tilde{f}, U_{\alpha}^{N}\tilde{g})_{F_N} &=& (U_{\alpha}^{N}\tilde{f}, \tilde{g})_{F_N}
\end{eqnarray*}
That is $U_{\alpha}^{\infty}$ is symmetric if all of the $U_{\alpha}^{N}$ are symmetric. Note that to get the last line in the calculation we used the fact that $\Phi^{*}_n$ is a measure preserving map from $\mathcal{B}(F_n)$ to $\mathcal{B}(L)$, which is a consequence of how the measures, $\mu_n$, are related to each other and to $\mu_{\infty}$.

Now it remains to show that from $U_{\alpha}^{0}$ being symmetric that $U_{\alpha}^{N}$ are all also symmetric. Already we have that  $U_{\alpha}^{\infty}$ being symmetric implies that $U_{\alpha}^{0}$ is symmetric. The symmetry of operators on collections of finite line segments is a well studied topic in Quantum Graph theory \cite{Kuchment2004, Kuchment2005, AkkermansEtAll2000}. By the way that the $U^{N}_{\alpha}$ were constructed inductively from $U^{0}_{\alpha}$ it is seen that all of the $U^{N}_{\alpha}$ are symmetric resolvents. 
\end{proof}

It is worth noting that the only facts that were used in proving this lemma were that we had a projective system of measure spaces, a family of resolvents satisfying Equation \ref{eqn:dynkin}, and facts about self-adjoint operators on quantum graphs. None of these things intrinsically are related to the production of a Laakso space and so this lemma is applicable in a much broader context than just this paper.

\begin{lemma}
If the sequence of spaces, $F_n$, are all quantum graphs and $X^{0}_t$ is  Feller process then $X_t$ is a Feller process.
\end{lemma}

\begin{proof}
That $X^{n}_t$ is Feller follows from the assumptions made on $X^{0}_t$ as part of the Barlow-Evans construction and from the fact that the $F_n$ are all quantum graphs. Since $X^{n}_t$ are Feller processes $U_{\lambda}^{n}:C(F_n) \rightarrow C(F_n)$. Now we use Equation \ref{eqn:dynkin} to say that $U^{\infty}_{\lambda} \Phi_n^{*} f \in C(L)$ for all $f \in C(F_n)$ for any $n \ge 0$. But in $C(L)$ the functions $\Phi_n^{*}C(F_n)$ are a dense subset so by taking uniform limits, because the resolvents are Markov, we get that $U^{\infty}_{\lambda}C(L) \subset C(L)$.
\end{proof}

\begin{lemma}\label{lem:feller}
If a process is Feller, then the associated Dirichlet form is regular.
\end{lemma}

This is Lemma 2.8 from \cite{Barlowetal2008}. We note this fact because we will be defining a Dirichlet form at the beginning of Section \ref{SharedProc} and proceed to show that certain continuous functions are dense in its domain. If the Dirichlet form were not already known to be regular the argument would be more delicate. 

It will be useful to fix some notation for function spaces that will be used to describe the domains of the operator and Dirichlet form associated to the process $X_t$. 

\begin{defn}\label{BF-FuncSpac} Recall that $\Phi_n:L \rightarrow F_n$ is the projection from the space $L$ to the $n$th level quantum graph from the Barlow Evans construction, and $\Phi^{*}_n$ is the pull back operator sending a function on $F_n$ to a function on $L$.

\begin{enumerate}
	\item Let $\tilde{A}_0 = \frac{\partial^{2}}{\partial x^{2}}$ on $[0,1]$ with Neumann boundary conditions.
	\item Let $\tilde{A}_n$ be the infinitesimal generator of the resolvent $U^{n}_{\alpha}$ that is associated to the process $X^{n}_t$ on $F_n$, which is by the construction in \cite{BarlowEvans2004} and \cite{EvansSowers2003} has the same action as $\tilde{A}_0$ on each line segment. Denote by $Dom(\tilde{A}_n) \subset L^{2}(F_n,\mu_n)$ domain of $\tilde{A}_n$.
	\item Let $G_n \subset C(F_n)$ be functions on the quantum graph, $F_n$, that are twice differentiable on each line segment of the graph, have continuous first and second derivatives on each line segment, and satisfy the Kirchoff matching conditions at each vertex. The Kirchoff matching condition states that the directional first derivatives along all the line segments meeting at a vertex sum to zero, see \cite{Kuchment2004}. The continuity condition implies that elements of $G_n$ are bounded as are their first and second derivatives over all of $F_n$.
	\item Let $\mathcal{D}_n = \Phi^{*}_n Dom(\tilde{A}_n) \subset L^{2}(L,\mu_{\infty})$ be the pull back of the domain of $\tilde{A}_n$ to a function space on $L$. This is so that the domains are subspaces of the same $L^{2}$ space.
	\item Let $\tilde{\mathcal{G}}_n$ be the set of  continuous functions on $F_n$ that are continuously differentiable on each line segment in $F_n$ and the derivatives have finite limits at the vertices. Set $\tilde{\mathcal{G}} = \bigcup_{n=0}^{\infty} \Phi_n^{*} \tilde{\mathcal{G}}_n$.
\end{enumerate}
\end{defn}

{\bf Remark:} \label{MPtoDF} There is a non-trivial Dirichlet form on the fractal $L$. Since there is a non-trivial symmetric Markov process, namely standard Brownian motion, which can be used in the Barlow-Evans construction there is a non-trivial symmetric Markov process on the fractal $L$ by the previous lemma. Which by Theorem \ref{DFMP-Cor2}\ yields a Dirichlet form which will be generated by a non-trivial self-adjoint linear operator, $\tilde{A}$. This operator is the generator of the resolvent $U^{\infty}_{\alpha}$. 

{\bf Remark:} The spaces $\C{G}$ and $\tilde{\C{G}}$ are the same function space on $L$. This is easily seem by tracing Definition \ref{BF-FuncSpac} through the homeomorphism, $\eta$, to Definition \ref{def:calG}. It is straight forward to see that the definitions are equivalent.

In Proposition \ref{1D-BR} and the last remark we know that there is a Laplacian on $L$ defined as the infinitesimal generator of the Markov process through Barlow and Evans' construction, i.e. $-\tilde{A}$. We also know that negative second differentiation with Neumann boundary conditions is the operator associated to one dimensional reflecting Brownian motion on the unit interval and that this is reminiscent of the operator $-A$ defined in the remark on \pageref{rem:A} which generates the Dirichlet form for the Laakso construction using minimal generalized upper gradients. To begin the process of showing that these are the same operators we look into the domain of $\tilde{A}$ with the intention of showing it is the same as the domain of $A$.

We will shortly be considering the closure of function spaces and of self-adjoint operators. For both of these the graph norm $\| u \|_{L^{2}(L,\mu_{\infty})} +  \| \tilde{A}u \|_{L^{2}(L,\mu_{\infty})}$ gives the relevant topology.

\begin{prop}\label{prop:AtFacts} Let $\tilde{A}_n$ as above, then
\begin{enumerate}
	\item For $f \in G_n$ $\tilde{A}_n f = \frac{\partial^{2}}{\partial x^{2}}f$ where $\frac{\partial^{2}}{\partial x^{2}}f$ is the second derivative of $f$ at each point in the interior of the line segments and not defined at the vertices which are a set of measure zero, and $Dom(\tilde{A}_n) = \overline{G_n}$ (closure taken in the graph norm),
	\item For $n \ge 0$, $\Phi^{*}_n G_n \subset \Phi^{*}_{n+1} G_{n+1}$,
	\item For $n \ge 0$ and $f \in G_n$,  $\tilde{A}\Phi^{*}_n f = \Phi^{*}_n \tilde{A}_n f$.
\end{enumerate}
\end{prop}

\begin{proof} We take each claim separately.
\begin{enumerate}
	\item By Definition \ref{BF-FuncSpac} part 2, $\tilde{A}_n$ acts on each line segment in $F_n$ in the same manner as $\tilde{A}_0$, which is the standard Laplacian on the line. For $f \in G_n$ when restricted to a line segment in $F_n$ is in the domain of the standard Laplacian and mapped to $\frac{\partial^{2}}{\partial x^{2}}f$ restricted to that line segment. The self-adjointness is given by the general theory in \cite{Kuchment2004} but can also be seen by using integration by parts on each line segment in $F_n$ and using the matching conditions built into the definition of $G_n$ to make the boundary terms vanish. So taken together we get the claim on all of $F_n$. 
	\item It is sufficient to show that $\phi^{*}_{n+1,n} f \in G_{n+1}$ for $f \in G_n$. As a pullback through a map as constructed in Section \ref{BEConst} it can been seen that $\phi^{*}_{n+1,n}f$ meets the criteria for membership in $G_n$. 
	\item This follows from the defining relationship of $U^{\infty}_{\alpha}$ that was given in Lemma \ref{lemma:symmetry} and the strong continuity of the resolvent of a Dirichlet form that makes it possible to relate $U^{\infty}_{\alpha}$ to its generator by a limit in the strong topology. 
\end{enumerate}
\end{proof}

\begin{prop}\label{prop:resolvent}
Let $U^{n}_{\lambda}$ be the resolvent associated to $\tilde{A}_n$ then $U^{n}_{\lambda}(\tilde{\mathcal{G}}_n) \subset \tilde{\mathcal{G}}_n$. 
\end{prop}

\begin{proof}
Since $X^{n}_t$ is a Feller process we already have that $U^{n}_{\lambda}(\tilde{\mathcal{G}}_n) \subset C(F_n)$. Let us pick $f \in \tilde{\mathcal{G}}_n$. Then $U^{n}_{\lambda}f = g$ for some $g \in Dom(\tilde{A}_n)$. But this is the same as saying that $f = \tilde{A}_ng+ \lambda g$. Due to Proposition \ref{prop:AtFacts} we know that $\tilde{A}_n$ is a local operator that on each line segment of the graph $F_n$ acts as second differentiation. With this locality we can look to see what properties $g$ has on each line segment individually. For $g$ to be an element of $\tilde{\mathcal{G}}_n$ is has to be continuous on $F_n$ which we already have and it also has to be continuously differentiable on each line segment with finite limits at the ends of the line segments. Compare with the comments in the proof of Theorem 17 in \cite{Kuchment2004}. The question now is if $f \in C^{1}([a,b])$ with bounded derivative the questions is if $g$ is as well when we have the relationship $f = g''+\lambda g$. But this is a standard question in the theory of ordinary differential equations and is known to be true. Thus $g \in \tilde{\mathcal{G}}_n$.

\end{proof}

The Laplacian $\tilde{A}$ is defined as the projective limit of the operators $\tilde{A}_n$, as a consequence of this definition $$\bigcup_{n=0}^{\infty} \mathcal{D}_n \subset Dom(\tilde{A})$$ is a dense subset in the graph norm $\| u \|_{L^{2}(L,\mu_{\infty})} +  \| \tilde{A}u \|_{L^{2}(L,\mu_{\infty})}$. We distinguish the projective limit of $\tilde{A}_n$ from the closure of an increasing family of self-adjoint operators because in themselves the domains of $\tilde{A}_n$ are in $C(F_n)$, see \cite{Lax2002} for closing self-adjoint operators. We call $\tilde{A}$ the \emph{projective limit} of $\tilde{A}_n$ if $(\tilde{A}, \mathcal{D}_n) = (\tilde{A}_n, \mathcal{D}_n)$ for all $n \ge 0$ and the operator $(\tilde{A}, \overline{\cup_n \mathcal{D}_n})$ is self-adjoint.

\begin{theorem}\label{thm:domtildea}
Using the graph norm $\| u \|_{L^{2}(L,\mu_{\infty})} +  \| \tilde{A}u \|_{L^{2}(L,\mu_{\infty})}$ on the space $Dom(\tilde{A})$ to define a topology,  we have that $\Phi^{*}_nG_n \subset \mathcal{D}_n \subset Dom(\tilde{A}),$
and
$$Dom(\tilde{A}) = \overline{ \bigcup_{n=0}^{\infty} \Phi^{*}_nG_n}.$$
\end{theorem}

\begin{proof}
The first claim is $\Phi^{*}_nG_n \subset \mathcal{D}_n \subset Dom(\tilde{A})$, which reduces to showing $G_n \subset Dom(\tilde{A}_n)$. On $F_n$ the boundary consists of vertices of degree one so the Kirchoff matching condition that elements of $G_n$ satisfy force the directional derivatives at all boundary points to be zero, so elements of $G_n$ satisfy the boundary conditions of $Dom(\tilde{A}_n)$. The action of $\tilde{A}_n$ is second differentiation on each line segment, the vertices being a null set can be set aside,   so for $f \in G_n$ the function $\frac{\partial^{2}f}{\partial x^{2}}$ where $x$ is a coordinate in any of the line segments is well defined and in $L^{2}(F_n,\mu_n)$ by the boundedness of the second derivatives imposed by the definition of $G_n$.

The second claim $\mathcal{D}_n = \overline{\Phi^{*}_nG_n}$, which reduced to $Dom(\tilde{A}_n) = \overline{G_n}$ is from the first part of Proposition \ref{prop:AtFacts}. The last claim that $Dom(\tilde{A}) =  \overline{ \bigcup_{n=0}^{\infty} \Phi^{*}_nG_n}$ holds because the associated Dirichlet form is regular by Lemma \ref{lem:feller}. 
\end{proof}

\section{A Shared Markov Process}
\label{SharedProc}
We have shown that the Laakso construction of $L$ guarantees that there is a Dirichlet form linked to the minimal generalized upper gradients which corresponds to some Markov process. We have seen from the Barlow-Evans construction that a Markov process is guaranteed to exist as well. The choice of base process as reflected Brownian motion was not the only possible decision. Other processes on $L$ could be built from Markov processes on the base space $F_0$ satisfying the assumption at the beginning on Section \ref{ProcBE}. These Markov processes give rise to generators which then give rise to Dirichlet forms. But as we have chosen a particular one let us stay with it and complete the comparison between the process $X^{\infty}_t$ and the Dirichlet form $(\C{E},Dom(\C{E}))$.

\begin{defn}\label{def:tildeE}
Let $\tilde{\mathcal{E}}$ be the Dirichlet form associated to the Markov process considered in Section \ref{ProcBE} via the self-adjoint operator, $\tilde{A}$, by the formula $$\tilde{\mathcal{E}}(u,u) = -\int_L u (\tilde{A}u)\ d\mu$$
For $u \in Dom(\tilde{A})$. The domain of $\tilde{\mathcal{E}}$ is $Dom(\sqrt{-\tilde{A}})$ as defined by the functional calculus for self-adjoint operators. 
\end{defn}

This way of associating a Dirichlet form and operator is the same as in Theorems \ref{SFSG-Cor}\ and \ref{DFMP-Cor}. It is worth noting that since the self-adjoint operator in question is the infinitesimal generator of a Markov process it is Markovian itself hence $\tilde{\mathcal{E}}$ is a Dirichlet form by Theorems \ref{DFMP-Cor}\ and \ref{SGMP-cor}. Recall the Dirichlet form $\mathcal{E}$ defined in Lemma \ref{LaaDF}, these two Dirichlet forms have their domains contained in $L^{2}(L)$ and to check that whether they are the same we have to first check that their domains have a common dense subset and then that they agree on this dense subset.

Before continuing to show that $\mathcal{E}$ and $\tilde{\mathcal{E}}$ are the same Dirichlet form we need a lemma to describe the domain of $\tilde{\mathcal{E}}$ in a manner that will be comparable to the description of the domain of $\mathcal{E}$ using the functions in $\mathcal{G}$ in Lemma \ref{LaaDF}.

\begin{lemma}
Let \C{G} is a dense subset of $Dom(\C{G})$ in the norm $\|u\|_{L^{2}} + \tilde{\mathcal{E}}(u,u)$.
\end{lemma}

\begin{proof}
With a similar argument as in Lemma \ref{StoneWeierstrass} we can see that \C{G} is a dense subset of $C(L)$ in the uniform norm, but $L$ is a finite measure space so $C(L) \subset L^{2}(L)$ and by standard results is a dense subset in the $L^{2}$ norm. Because of this $U^{\infty}_{\lambda}(\C{G})$ is a dense subset of $Dom(\tilde{A})$ in it's topology because the resolvents are continuous maps. It is also known that $Dom(\tilde{A})$ embeds continuously as a dense subset of $Dom(\sqrt{-\tilde{A}})$. This leaves only whether or not $U^{\infty}_{\lambda}(\C{G}) \subset \C{G}$ but this has been shown in Proposition \ref{prop:resolvent}.
\end{proof}

Now we are ready to state and prove the main result of the paper.

\begin{theorem}\label{thm:DomEq}
The two Dirichlet forms, $\mathcal{E}$ and $\tilde{\mathcal{E}}$, are equal.
\end{theorem}

\begin{proof}
In order to show that the two Dirichlet forms, $\mathcal{E}$ and $\tilde{\mathcal{E}}$, are equal we show that they agree on a dense subset of their domains. Since we already know this subset is dense in the domains of $\mathcal{E}$ and $\tilde{\mathcal{E}}$ this agreement will extend to their full domains by the \emph{same} metric. Hence the two Dirichlet forms will have the same domains and give the same values to functions in their common domain. 

The subset of $Dom(\mathcal{E})$ and $Dom(\tilde{\mathcal{E}})$ that we will consider is \C{G}. As remarked after Definition \ref{BF-FuncSpac}, $\C{G} = \tilde{\C{G}}$ is dense in both domains. Both $\mathcal{E}$ and $\tilde{\mathcal{E}}$ are the integrals of first derivatives squared on functions in \C{G}. By Theorem \ref{thm:diffu} we give meaning to this in language of Laakso's construction. At the beginning of this Section we defined $\tilde{\mathcal{E}}$ in terms of the self-adjoint operator $\tilde{A}$ and said that it's domain is $Dom\left(\sqrt{-\tilde{A}}\right)$ without saying what was in that domain. However using the results in \cite{Kuchment2004} we see that $\sqrt{-\tilde{A}_n}$ on $F_n$ with domain $\overline{\mathcal{G}_n}$ with the closure taken in the metric given by 
$$\| \cdot \|_2 + \left\| \sqrt{-\tilde{A}_n} \cdot \right\|_2 = \| \cdot \|_2 + \tilde{\mathcal{E}}(\cdot, \cdot).$$

By Definitions \ref{def:calG} and  \ref{BF-FuncSpac}.5 and the definitions of $\mathcal{E}$ and $\tilde{\mathcal{E}}$ we can see that for $f \in \mathcal{G}$ that $\mathcal{E}(f,f) = \tilde{\mathcal{E}}(f,f)$. Since the two Dirichlet forms have the same domains and agree on a dense subset they are the same.
\end{proof}

\begin{cor} The Dirichlet form from Theorem\ \ref{LaaDF} associated to the minimal generalized upper gradients on the fractal $L$ corresponds to the Markov process from the Barlow-Evans construction with $X^{0}_t$ the standard Brownian motion on the unit interval.
\end{cor}

\begin{proof}
By Theorem \ref{thm:DomEq}\ the Dirichlet forms generated by the minimal generalized upper gradients and to the Markov processes built from Brownian motion on the unit interval are associated to the same Dirichlet form. Then by Theorem \ref{DFMP-Cor} the Dirichlet forms, Markov processes, and self-adjoint operators from both constructions are the same, $A = \tilde{A}$.
\end{proof}

{\bf Remark:} It is unlikely that this sort of result would hold for a general space constructed by Barlow and Evan's construction. What makes it possible in this situation is the well defined cell structure where the interior and complement of a cell are disjoint. As well as having approximating 1-dimensional spaces. Which keeps the possible sorts of processes relatively accessible objects with which to work. It seems reasonable that as long as a cell structure is available that this type of method should be possible for more families of Barlow-Evans spaces.

\begin{center}{\small
Contact:\\
steinhurst@math.uconn.edu\\
Department of Mathematics\\
University of Connecticut\\
Storrs~CT~06269~USA}\end{center}

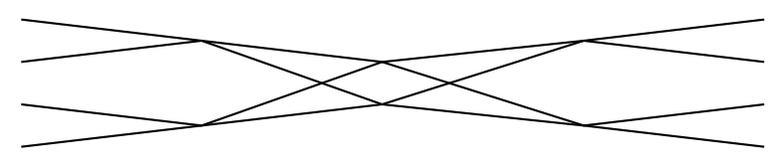
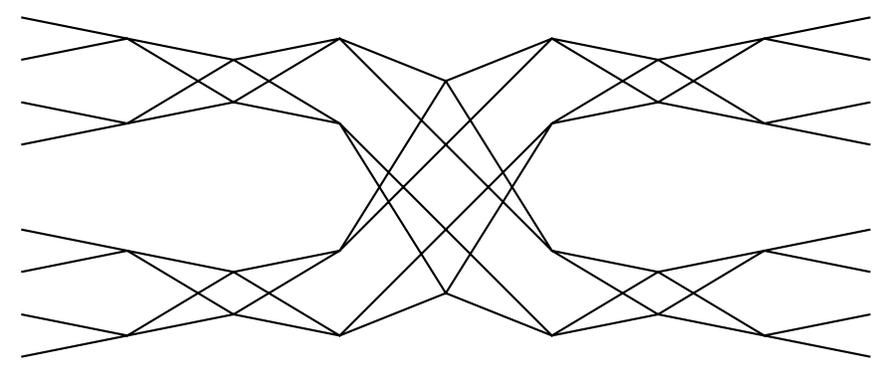
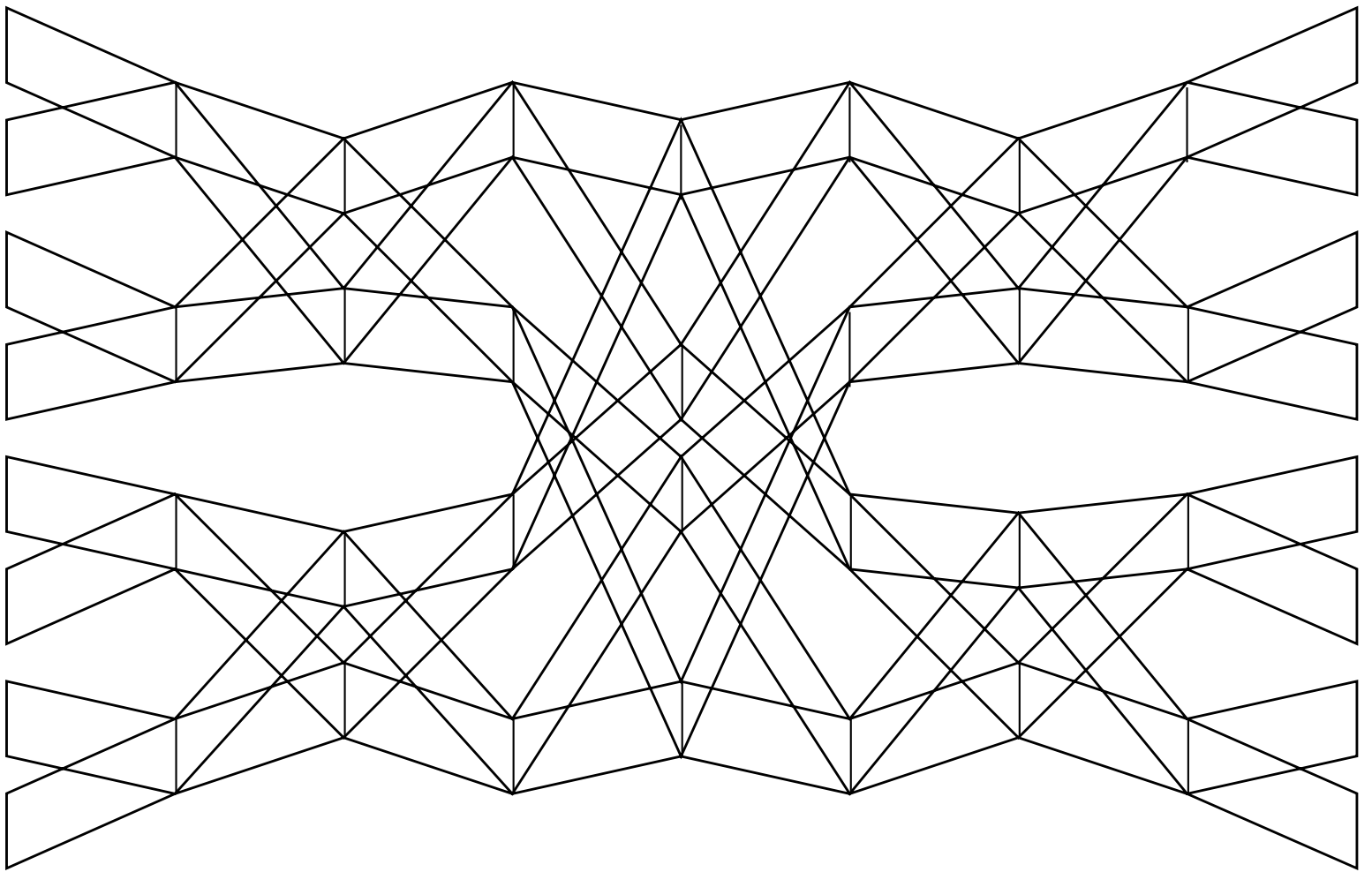
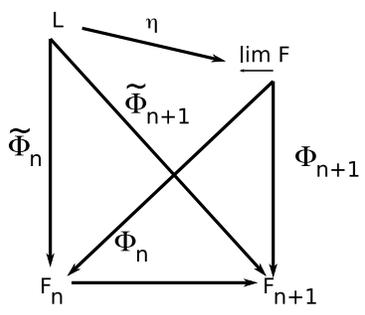
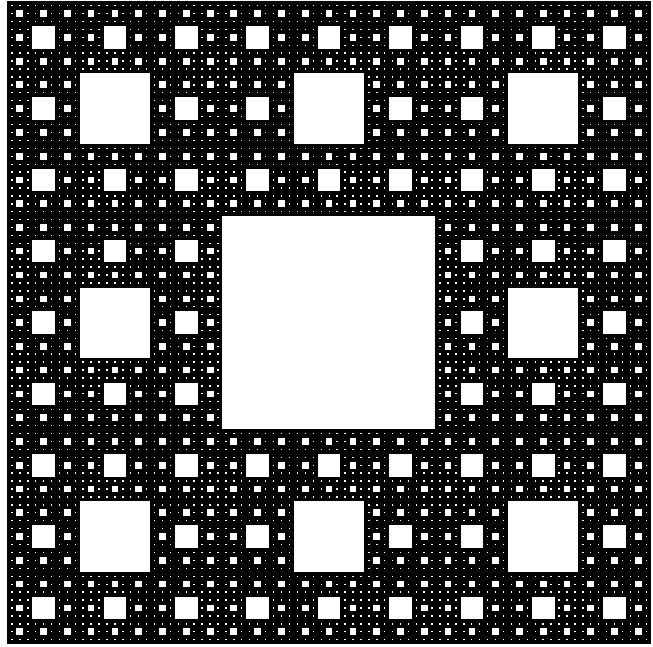
\end{document}